\documentclass[a4paper,11pt]{article}
\usepackage{graphicx}
\usepackage{amssymb}
\usepackage{amsmath}
\usepackage{booktabs}
\usepackage{subfig}
\usepackage{lmodern}
\usepackage{natbib}
\usepackage{tikz}

\newtheorem{theorem}{Theorem}
\newtheorem{definition}[theorem]{Definition}

\newcommand{\degree}{\ensuremath{^\circ}}
\newcommand{\field}[1]{\mathbb{#1}}  
\newcommand{\R}{\field{R}} 
\newcommand{\Zset}{\mathcal{Z}} 
\newcommand{\transp}{^{^{\top}}}
\newcommand{\twonorm}[1]{\left|\left|#1\right|\right|_2}
\newcommand{\infnorm}[1]{\left|\left|#1\right|\right|_\infty}
\renewcommand{\vec}[1]{\boldsymbol{#1}} 

\newcommand{\ea}[1]{\ell_{#1}(\alpha)}

\newcommand{\IPS}{IPS}
\newcommand{\OPS}{OPS}
\newcommand{\PS}{PS}

\hypersetup{
    colorlinks=true,
    citecolor=black,
    urlcolor=blue,
    linkcolor=black,
    pdfborder={0 0 0},
}

\begin{document}

\title{Approximating the Pareto set of multiobjective linear programs via Robust Optimization}
\author{Bram L. Gorissen, Dick den Hertog \\ \\
\textit{\small Tilburg University, Department of Econometrics and Operations Research, 5000 LE Tilburg, Netherlands} \\
\textit{\small {\tt \{b.l.gorissen,d.denhertog\}@tilburguniversity.edu}}}
\date{}
\maketitle

\begin{abstract}
We consider problems with multiple linear objectives and linear constraints and use Adjustable Robust Optimization and Polynomial Optimization as tools to approximate the Pareto set with polynomials of arbitrarily large degree. The main difference with existing techniques is that we optimize a single (extended) optimization problem that provides a polynomial approximation whereas existing methods iteratively construct a piecewise linear approximation. One of the advantages of the proposed method is that it is more useful for visualizing the Pareto set.
\end{abstract}

\begin{tikzpicture}[remember picture,overlay]
\node[anchor=south,yshift=10pt] at (current page.south) {\fbox{\parbox{\dimexpr\textwidth-\fboxsep-\fboxrule\relax}{\footnotesize This is an author-created, un-copyedited version of an article published in Operations Research Letters \href{http://dx.doi.org/10.1016/j.orl.2012.05.007}{DOI:10.1016/j.orl.2012.05.007}.}}};
\end{tikzpicture}

\section{Introduction}
Multiobjective optimization problems (MOPs) have a vector valued objective function $\vec{f} = [f_i]_{1 \leq i \leq k}$, where each $f_i$ is a separate objective. Often it is not possible to have optimal values for all $f_i$ simultaneously, e.g. in portfolio optimization it is not possible to have minimum risk and maximum return at the same time. Another example is intensity-modulated radiation therapy, where tumour coverage is balanced with sparing of surrounding organs \citep{Craft2008,RennenSW}. Optimization of a vector valued function involves a trade-off between two or more objectives $f_i$ ($1 \leq i \leq k$). 

A simple way to deal with multiple objectives is by assigning an importance factor $w_i >0$ to each objective and optimizing $\sum_{i=1}^k w_i f_i$ (we make the assumption that all $f_i$ should be minimized w.l.o.g.). If such importance factors are not known a priori, a Pareto set ($\PS$) allows the decision maker to make the trade-off after optimization. The set $\PS$ consists of all objective vectors $\vec{f}$ in which one or more objectives can not be improved without deteriorating one or more other objectives. Overviews of MOPs and approximation methods can be found in \citep{Branke2008,Ehrgott2005,miettinen1999}.

In practice often approximations of $\PS$ are used, since the exact $\PS$ can often not be found. In literature many different approximation methods are proposed. 
It is desirable to approximate $\PS$ with as few optimization runs as possible \citep{RennenSW}.
A well-known class of such approximation methods is sandwich methods.
Sandwich methods \citep{RennenSW} produce  piecewise linear approximations in which between $\PS$ is located. In each iteration an optimization problem is solved, which leads to adding one or more facets to the approximations.
All other methods (again see \citep{Branke2008,Ehrgott2005,miettinen1999}) are also sequential, i.e. in each iteration one has to solve an optimization problem which leads to an improvement of the approximation.

In this paper, we  focus on approximating $\PS$ for linear programming and propose a totally different way than those in the literature.
The first difference is that our method is not sequential, but generates the approximation by solving one extended optimization problem.
The second difference is that the final approximation is not piecewise linear but a polynomial.
The way we construct this approximation is by using techniques from Adjustable Robust Optimization (ARO) \citep{BenTal} and Polynomial Optimization \citep{Laurent2009}. We first explain the link to ARO.
The Pareto set is seen as a function of the {\it uncertain parameters}  $f_1,..., f_{k-1}$. The area of interest, i.e. the domain for $f_1,..., f_{k-1}$ for which we would like to approximate the Pareto set, is considered as the 
{\it uncertainty region}. All variables in the linear program are made {\it adjustable} in the parameters $f_1,..., f_{k-1}$. We use  polynomials for the decision rule, and use  Polynomial Optimization theory to reformulate the resulting robust counterpart into a Semi-Definite Programming (SDP) problem. Since the number of uncertain parameters (i.e. $k-1$) is often low, the sizes of the LMIs in the SDP are relatively small.
Notice that in our approach ARO is merely used as a tool, uncertainty in the data is not considered.

The approach proposed in this paper has the following advantages:

The first advantage is that the final approximation is more tractable for navigating through $\PS$. The polynomial representation is useful for the user to visualize $\PS$ for selecting the final solution.  This is the reason why in \citep{Goel2007}, afer determining points in the feasible region close to or on $\PS$, polynomial regression is used to obtain a tractable representation of $\PS$ (a so-called response surface). 
Our method finds such a tractable representation directly, with the additional advantage that it is guaranteed to lie in the feasible region.
Sometimes the decision maker needs a local approximation of $\PS$ around a given solution. \cite{zhang2000} formulate and test a method that gives a local quadratic approximation of (not necessarily convex) $\PS$, but this approximation is neither an inner nor an outer approximation. Our method gives a polynomial of arbitrary degree and is guaranteed to give an inner approximation.

The second advantage is that our approach can be used to certify with a single optimization run that a given set $V$ is dominated, i.e. that all elements of the set are dominated.
If our method finds a feasible solution, then this solution is a certificate that the set $V$ is dominated.

The third advantage is that the explicit polynomial approximation can be used in optimization problems with Pareto constraints. Such problems contain constraints that enforce that the solution should be (near) Pareto optimal for a certain multi-objective linear program. Examples of these problems can be found in e.g. \citep{Hackman2002}.
A Pareto constraint can be replaced by the explicit polynomial approximation found by our method.

The fourth advantage is the possibility to quickly determine the shape of $\PS$. 
\citep{Craft2008}, e.g., show that in IMRT the set of feasible objective vectors is often ``long and narrow'' and therefore a linear approximation of $\PS$ suffices.  This linear approximation can be easily obtained by our approach. 
Finally, after an initial approximation, the most interesting subregion can be selected, followed by one or more iterations of approximation and selection. An interesting subregion can also be used as input for another algorithm that explores it more carefully.

This method also has five disadvantages. First, the resulting problem is often an SDP while the original problem is LP. Only in case the approximation is linear and the region of interest is polyhedral, the resulting problem is LP. Note that this is still an interesting case; see \citep{Craft2008} that uses linear approximations in IMRT problems.  Second, our method requires the region of interest to be known. 
Sandwich algorithms are capable of exploring the region of interest \citep{RennenSW}. Third, it is difficult to approximate the Pareto set at its vertices, because polynomials are smooth functions. The further the vertex angle from 180\degree, the more difficult it is to approximate it. However, in some cases a smooth approximation is desirable, see e.g. \citep{Mello2002}. Fourth, the method can not be extended to nonlinear multiobjective problems with current ARO technology. Methods for approximating nonlinear MOPs can be found in \citep{Luque2012,Utyu2005}. Fifth, while we show (Appendix \ref{sec:outerapprox}) that the method can also produce outer approximations, this is practically impossible due to computational issues.

\section{Notation}
We use the notation from \citep{RennenSW} with some minor changes.

Throughout this paper, we use the following orderings of vectors. Let $\vec{x},\vec{y}\in\R^{n}$ with $n\geq2$. With $x_{i}$, we denote the $i^{th}$ element of the vector $\vec{x}$. To enumerate different vectors, we use superscripts. When ordering two vectors, we use:
\begin{itemize}
\item $\vec{x}<\vec{y} \Leftrightarrow x_{i}<y_{i}$ for all $i=1,\ldots,n$.
\item $\vec{x}\lneqq \vec{y} \Leftrightarrow x_{i}\leq y_{i}$ for all $i=1,\ldots,n$ and $\vec{x}\neq \vec{y}$.
\item $\vec{x}\leq \vec{y} \Leftrightarrow x_{i}\leq y_{i}$ for all $i=1,\ldots,n$.
\end{itemize}
The symbols $>$, $\gneqq$ and $\geq$ are defined accordingly. We furthermore define the set $\R^{n}_{-}= \{\vec{x}\in\R^{n}:\vec{x}\leq0\}$. If $X\subseteq\R^{n}$, then we define $X+\R^{n}_{-}=\{\vec{x}+\vec{y} : \vec{x} \in X, \vec{y} \in \R^n_-\}$. The sets $\R^{n}_{+}$ and $X+\R^{n}_{+}$ are defined accordingly.

In this paper, we consider the following multi-objective optimization problem:
\begin{align*}
 \min        \qquad & \vec{f(x)}=[(\vec{c^1}) \transp \vec{x},\ldots,(\vec{c^k}) \transp \vec{x}]^{\top} \\
                    & A\vec{x} \leq \vec{b},
\end{align*}
where $\vec{x} \in \R^n$ is the optimization variable, $\vec{c^i} \in \R^n$ are the objective vectors, and $A \in \R^{m \times n}$ and $\vec{b} \in \R^m$. 

As it is generally impossible to find a feasible $x$ that minimizes all objectives at the same time, our aim is to find a set of so-called Pareto optimal solutions.

\begin{definition}
$\!$\newline An objective vector $\vec{f(x)}$, for $\vec{x}$ such that $A\vec{x}\leq \vec{b}$, is (strongly) dominated if there exists an $\vec{\widetilde{x}}$ such that $A\vec{\widetilde{x}}\leq \vec{b}$ and $\vec{f(\widetilde{x})}< \vec{f(x)}$. If no such $\vec{\widetilde{x}}$ exists, the objective vector $\vec{f(x)}$ is weakly Pareto optimal.\\
An objective vector $\vec{f(x)}$, for $\vec{x}$ such that $A\vec{x} \leq \vec{b}$, is weakly dominated if there exists an $\vec{\widetilde{x}}$ such that $A\vec{\widetilde{x}}\leq \vec{b}$ and $\vec{f(\widetilde{x})}\lneqq \vec{f(x)}$. If no such $\vec{\widetilde{x}}$ exists, the objective vector $\vec{f(x)}$ is (strongly) Pareto optimal.
\end{definition}

The set of Pareto optimal solutions is denoted by $\PS$. An inner and outer approximation of the Pareto set are defined as follows:
\begin{definition}
A set $\IPS\subseteq \R^k$ is an inner approximation of $\PS$ if it satisfies $\IPS\subseteq \PS+\mathbb{R}^{k}_{+}$.
\end{definition}
\begin{definition}
A set $\OPS\subseteq \R^k$ is an outer approximation of $\PS$ if it satisfies $\PS\subseteq \OPS+\mathbb{R}^{k}_{+}$.
\end{definition}
We will approximate $\PS$ with polynomials on multidimensional sets. The following definitions are used to define the degree of a polynomial.
\begin{definition}
A monomial of degree $d$ in $\vec{x} \in \R^n$ with powers $\vec{a} \in \R^n$ such that $\sum_{i=1}^n a_i = d$, is defined by $\prod_{i=1}^n x_i^{a_i}$.
\end{definition}
\begin{definition}
A polynomial of degree $d$ in $\vec{x} \in \R^n$ is defined as the sum of monomials in $\vec{x}$ of degree up to $d$. The degree of a polynomial $f$ is denoted as $deg(f)$.
\end{definition}
\section{Inner approximation} \label{sec:ia}
Let $U \subseteq \R^{k-1}$ be the domain of interest for \\ $((\vec{c^1}) \transp \vec{x},(\vec{c^2}) \transp \vec{x},\ldots,(\vec{c^{k-1}}) \transp \vec{x})$. For a fixed $\vec{u}$ in $U$, the following optimization problem determines a single weakly Pareto optimal solution \cite[Thm. 3.2.1]{miettinen1999}:
\begin{align*}
 \min_x      \qquad & (\vec{\vec{c^k}}) \transp \vec{x} \\
                    & (\vec{c^i}) \transp \vec{x} \leq u_i && i=1,2,\ldots,k-1  \\
                    & A\vec{x} \leq \vec{b}.
\end{align*}

If the solution $x$ is unique, it is (strongly) Pareto optimal \cite[Thm. 3.2.4]{miettinen1999}. For every $\vec{u}$, there will be a different optimal $\vec{x}$. So, we want to solve for $\vec{x}$ as a function of $\vec{u}$. The constraints should hold for all $\vec{x(u)}$ for which $\vec{u}$ is in $U$, and the goal is e.g. to minimize the average objective:

\begin{subequations}
\begin{align}
 \min_{\vec{x(u)}} \qquad & \int_{U} (\vec{c^k}) \transp \vec{x(u)} du  && \label{ia-obj} \\
                    & (\vec{c^i}) \transp \vec{x(u)} \leq u_i           && \forall \vec{u} \in U, \quad i=1,2,\ldots,k-1 \label{ia-moving} \\
                    & A\vec{x(u)} \leq \vec{b}                          && \forall \vec{u} \in U. \label{ia-feas}
\end{align}
\end{subequations}
This is an ARO problem, where $\vec{u}$ is the uncertain parameter, $U$ is the uncertainty region, and $x$ is an adjustable variable \cite{BenTal}. It is difficult to optimize over functions, therefore ARO uses parameterized functions for adjustable variables. The adjustable variables then become expressions that are linear in the parameters. For instance, if we take a linear parameterization $\vec{x(u)} = \vec{\alpha^0} + \alpha^1 \vec{u}$, the parameters are $\vec{\alpha^0} \in \R^{n}$ and $\alpha^1 \in \R^{n \times (k-1)}$. After substituting $\vec{x(u)}$ in the problem \eqref{ia-obj}-\eqref{ia-feas}, an ARO problem with constraints that are linear in $\vec{\alpha^0}$ and $\alpha^1$ remains. In general, the tractability of \eqref{ia-obj}-\eqref{ia-feas} depends on the class of functions considered for $x$ and the set $U$. Given a solution to this optimization problem, an inner approximation is given by $\newline \{ ((\vec{c^1}) \transp \vec{x(u)}, (\vec{c^2}) \transp \vec{x(u)}, ..., (\vec{c^k}) \transp \vec{x(u)}) : \vec{u} \in U \}$. Constraint \eqref{ia-feas} ensures that $\vec{x(u)}$ is feasible. So, the resulting inner approximation indeed lies in $\PS+\R^k_+$. The objective \eqref{ia-obj} minimizes the volume under the inner approximation if \eqref{ia-moving} is tight for all $\vec{u}$ . Note that a constant function $\vec{x(u)}$ may be feasible for this optimization problem. The reason why this problem returns an interesting inner approximation and not just a constant function $x$ is because from the objective it follows that a smaller $(\vec{c^k}) \transp \vec{x(u)}$ is better, and a smaller $(\vec{c^k}) \transp \vec{x(u)}$ can only be obtained by increasing $(\vec{c^i}) \transp \vec{x(u)}$ for $i=1,2,\ldots,k-1$. This increase as a function of $\vec{u}$ is constrained by \eqref{ia-moving}.

A functional description of the inner approximation may be more tractable for a specific purpose. Therefore, we define:
\[        \IPS = \{ (\vec{u}, (\vec{c^k}) \transp \vec{x(u)}) : \vec{u} \in U \}.      \]
If constraint \eqref{ia-moving} is tight for all $\vec{u}$, $\IPS$ is the same as the inner approximation given before. Otherwise, the given inner approximation dominates $\IPS$.

The question arises for which regions $U$ and functions $x$ this formulation is tractable. When $U$ is polyhedral, for instance a box ($\{\vec{u} \in \R^{k-1} : \infnorm{\vec{u}} \leq 1 \}$), and $x$ is linear in $\vec{u}$, this problem can be formulated as an LP. The optimization problem minimizes the volume enclosed between the linear approximation, the Pareto curve and the boundary of $U$. In case of two objectives ($k=2$), it will find the line connecting the point on the Pareto curve where $(\vec{c^1}) \transp \vec{x} = -1$ with the point on the Pareto curve where $(\vec{c^1}) \transp \vec{x} = 1$. Hence, it finds two Pareto optimal points. This extends to larger $k$, where it finds a plane going through $k-1$ Pareto optimal points. In case $U$ is a box, the instance size also grows linearly in $k$. This result of linear growth is new, because determining a linear inner approximation over a box would require determining the Pareto optimal points at the exponentially growing number of extreme points of $U$.

The inner approximation becomes more interesting when $x$ is nonlinear in $\vec{u}$. When $k=2$, a tractable choice is given by polynomials: for polynomials of arbitrarily large degree, the problem can be formulated as an SDP \cite[Lemma 14.3.4]{BenTal}, for which polynomial time solvers are available. When $k>2$, the problem is tractable when $U$ is ellipsoidal and $x$ is quadratic in $\vec{u}$ \cite[Lemma 14.3.7]{BenTal}. The resulting problem is an SDP with $m+k-1$ variables matrices of size $k+1$ and $m+k-1$ constraints.

For $k>2$, the problem can also be reformulated as an SDP when $U$ is a semialgebraic set ($U = \{ \vec{u} : p_i(\vec{u}) \leq 0 \quad (i \in I) \}$, where $p_i$ are polynomials of arbitrary degree), and $\vec{x}$ is a polynomial in $\vec{u}$, but the reformulation is not always equivalent. This means that the resulting optimal solution for the SDP reformulation may not be optimal for \eqref{ia-obj}-\eqref{ia-feas}. However, the solution to the reformulation is always an inner approximation, and the numerical results in Section \ref{sec:threeobj} are promising. The reformulation is based on polynomials that are sums of squares (SOS) of other polynomials. An example of an SOS polynomial is $5x^2 + 2x + 65$, because it can be decomposed as $(2x+4)^2 + (x-7)^2$. Testing whether a polynomial is SOS is equivalent to solving an SDP \cite[Lemma 3.8]{Laurent2009} with a matrix of size ${ k-1+d \choose d }$, where $k-1$ is the number of variables and $d$ is the degree of the polynomial. An SOS polynomial is obviously nonnegative. Let us focus on constraint \eqref{ia-moving}:
\[ u_1 - (\vec{c^1}) \transp \vec{x(u)} \geq 0 \qquad \forall \vec{u} : p_i(\vec{u}) \leq 0 \quad (i \in I). \]
By applying Putinar's Positivstellensatz \cite[Thm. 3.20]{Laurent2009}, which has been done before in ARO \cite{BertsimasPolynomials}, we can obtain a sufficient condition under which this constraint holds:
\begin{align}
u_1 - (\vec{c^1}) \transp \vec{x(u)} = \sigma_0(\vec{u}) + \sum_{i \in I}p_i(\vec{u}) \sigma_i(\vec{u}), \label{ia-moving-sos}
\end{align}
where $\sigma_0$ and $\sigma_i$ ($i \in I$) are SOS. Solving a problem with constraint \eqref{ia-moving-sos} instead of \eqref{ia-moving} is conservative for two reasons. First, \eqref{ia-moving-sos} may not be a necessary condition when there are no $\sigma_0, \sigma_i$ that are SOS for which the set $\{\vec{u} : \sigma_0(\vec{u}) + \sum_{i \in I}p_i(\vec{u}) \sigma_i(\vec{u}) \geq 0 \}$ is compact, or when \eqref{ia-moving} is a not a strict inequality. Compactness can easily be guaranteed by including a constraint $\sum_{i=1}^{k-1}u_i^2 - R \leq 0$ to the description of $U$, which can be done without changing $U$ because $U$ is bounded \cite[p. 186]{Laurent2009}. However, \eqref{ia-moving} will in general not be a strict inequality. Second, solving a problem with constraint \eqref{ia-moving-sos} as an SDP requires bounding the degree of $\sigma_0$ and $\sigma_i$.

We let the degree of $x$ determine the complexity of the problem unless $g$ is of higher degree, so we take the degree of $\sigma_0$ equal to $\max\{ deg(x), \max_i\{ deg(g_i) \} \}$, and the degree of $\sigma_i$ equal to $\max\{ 0, deg(\sigma_0) - deg(g_i) \}$.

An overview of all tractable cases is given in Table \ref{tbl:tractability}.

\begin{table}
\caption{Tractability of the inner approximation.}
\label{tbl:tractability}
\begin{center}
\begin{tabular}{rllll}
\toprule 
$k$          & $U$           & $\vec{x(u)}$ & Tractability & Exact \\
\midrule 
$\geq 2$     & box           & linear       & LP           & $\surd$ \\
$\geq 2$     & polyhedral    & linear       & LP           & $\surd$ \\
$\geq 2$     & ball          & linear       & CQP          & $\surd$ \\
$     2$     & interval      & polynomial   & SDP          & $\surd$ \\
$   > 2$     & ellipsoidal   & quadratic    & SDP          & $\surd$ \\
$\geq 2$     & semialgebraic & polynomial   & SDP          & $-$ \\
\bottomrule 
\end{tabular}
\end{center}
\end{table}
For many uncertainty regions it may be difficult to reformulate the integral in the objective function as a simple linear function in the optimization variables. In that case the objective can be replaced with the average value of $(\vec{c^k}) \transp \vec{x(u)}$ at well distributed sampling points $\vec{u}$ in $U$. For efficient sampling from a polytope, see e.g. \cite{Kannan2009}.

The user has to specify the domain of interest $U$. If the specified region is too large, two things may occur. First, $U$ may contain a vector with objective values that are too optimistic in the sense than they can not be met, in which case constraint \eqref{ia-moving} is infeasible. Second, $U$ may contain objective values that are not weakly dominated by $\vec{f(x)}$ for any feasible $\vec{x}$. In that case constraint \eqref{ia-moving} will not be tight and also the objective is not fully related to the area of interest.

When the number of objectives is three or more, a weak parameterization of $x$ can be another cause of infeasibility. If feasible solutions exist for all $\vec{u}$ in $U$, it is possible that these solutions can not be attained with the parameterization. An easy example is the case where $U$ contains two different Pareto optimal solutions (projected on the first $k-1$ coordinates) while $x$ is a constant function. When the optimization problem is infeasible, we do not see a possibility to detect whether the parameterization is too weak or $U$ contains infeasible points.

It is known that $\PS$ is convex, and that it is nonincreasing. It may be the case that these properties do not hold for the inner approximation, which is problematic when the inner approximation is used in an algorithm that assumes these properties to hold. In case $x$ is linear in $\vec{u}$, these conditions are automatically satisfied. In case $k=2$, $u$ is one-dimensional and convexity and nonincreasingness of $(\vec{c^2}) \transp \vec{x(u)}$ can be enforced by constraining the first and second derivative w.r.t. $u$. The first and second derivative of a polynomial is again a polynomial, so constraining these to be negative and positive, respectively, for polynomial $x$ does not increase the complexity class of the problem. In case $k>2$, $U$ is a polynomial and $x$ is quadratic in $\vec{u}$, say $x(\vec{u})_i = \alpha^0_i + (\vec{\alpha^1_i}) \transp \vec{u} + \vec{u} \transp \Gamma_i \vec{u}$ where $\alpha^0_i \in \R, \vec{\alpha^1_i} \in \R^{k-1}$ and $\Gamma_i \in \R^{(k-1)\times(k-1)}$ are decision variables. Nonincreasingness can easily be enforced by adding the constraints $\vec{\alpha^1_i} + 2 \Gamma_i \vec{u} \leq 0$ for all $\vec{u}$ in $U$, which is a set of $n(k-1)$ linear constraints with ellipsoidal uncertainty, each of which can be reformulated as a conic quadratic constraint. For convexity it is required that $\Gamma_i$ is positive semidefinite, which is an SDP constraint.

\section{Numerical examples}
\subsection{Two objectives} \label{sec:twoobj}
We construct a semi-random $150 \times 170$ matrix $A$, a 150-vector $b$, and two 170-vectors $\vec{c^1}$ and $\vec{c^2}$, such that the Pareto cuve is interesting on the interval $[0,25]$. We compute a polynomial inner approximations of degree up to 16.

\begin{figure}
	\centering
	\subfloat[degree 1]{
		\includegraphics[scale=0.35]{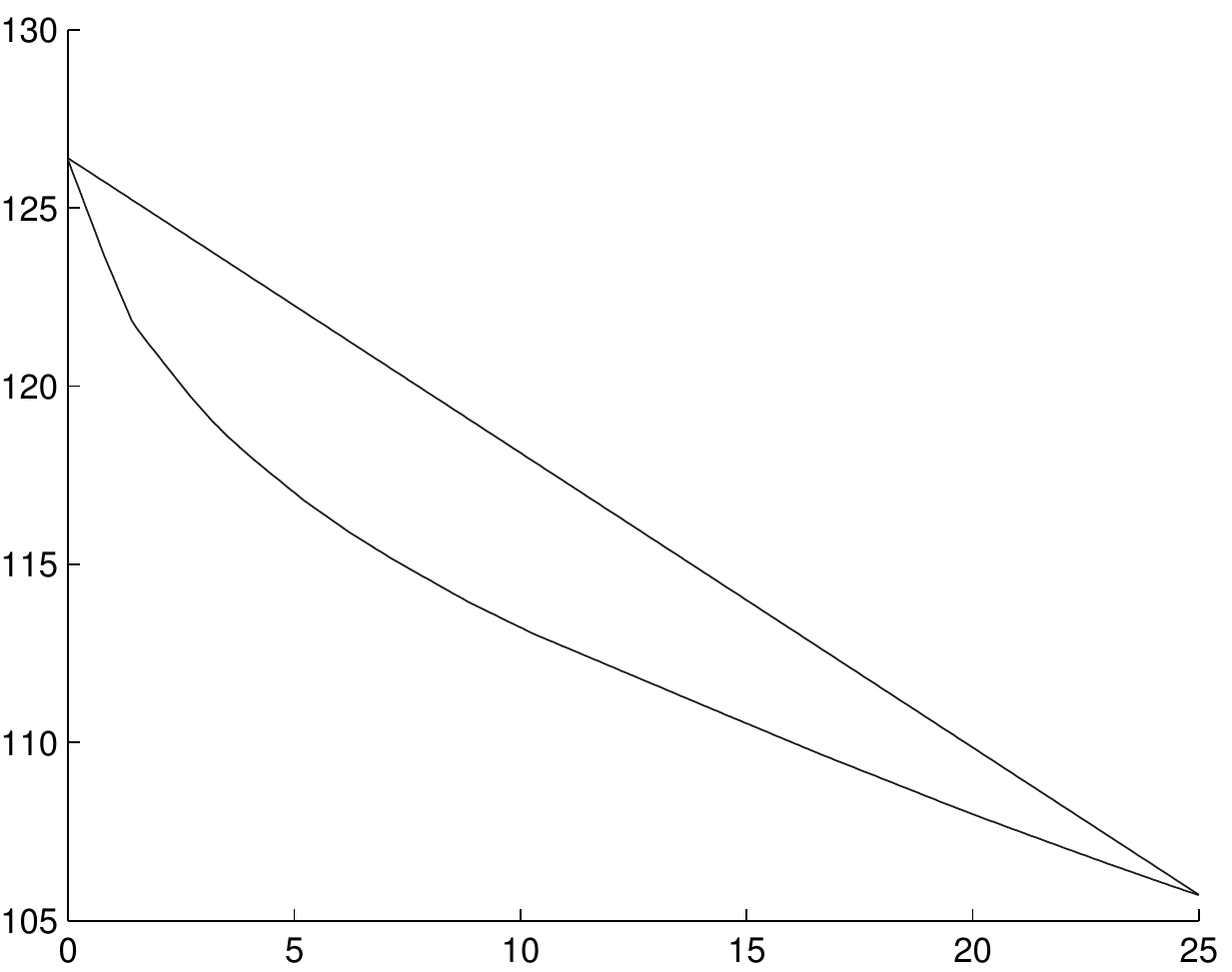}
		\label{fig:ex-dim2-1}
  }
	\subfloat[degree 4]{
		\includegraphics[scale=0.35]{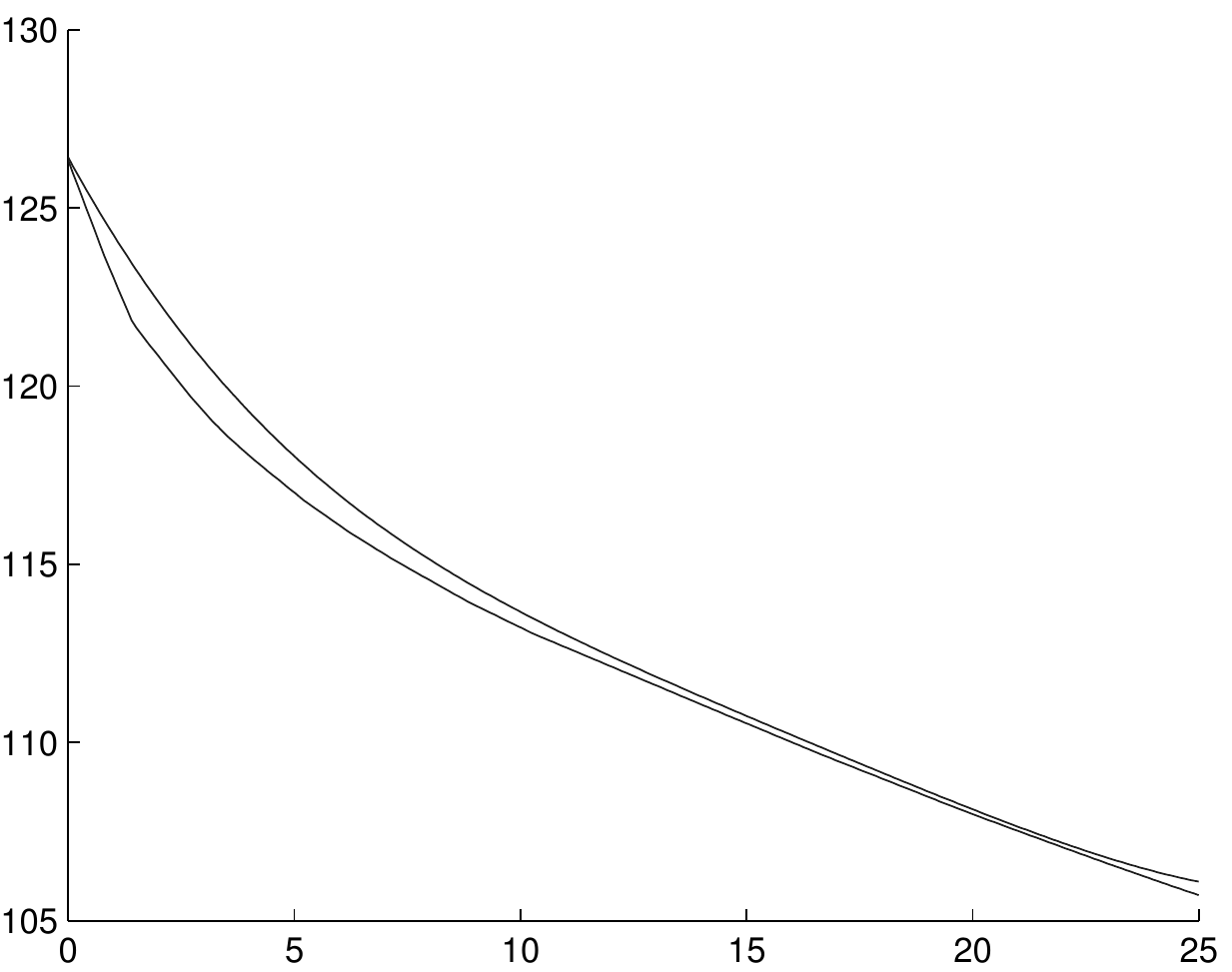}
		\label{fig:ex-dim2-2}
	}
	\\
	\subfloat[degree 8]{
		\includegraphics[scale=0.35]{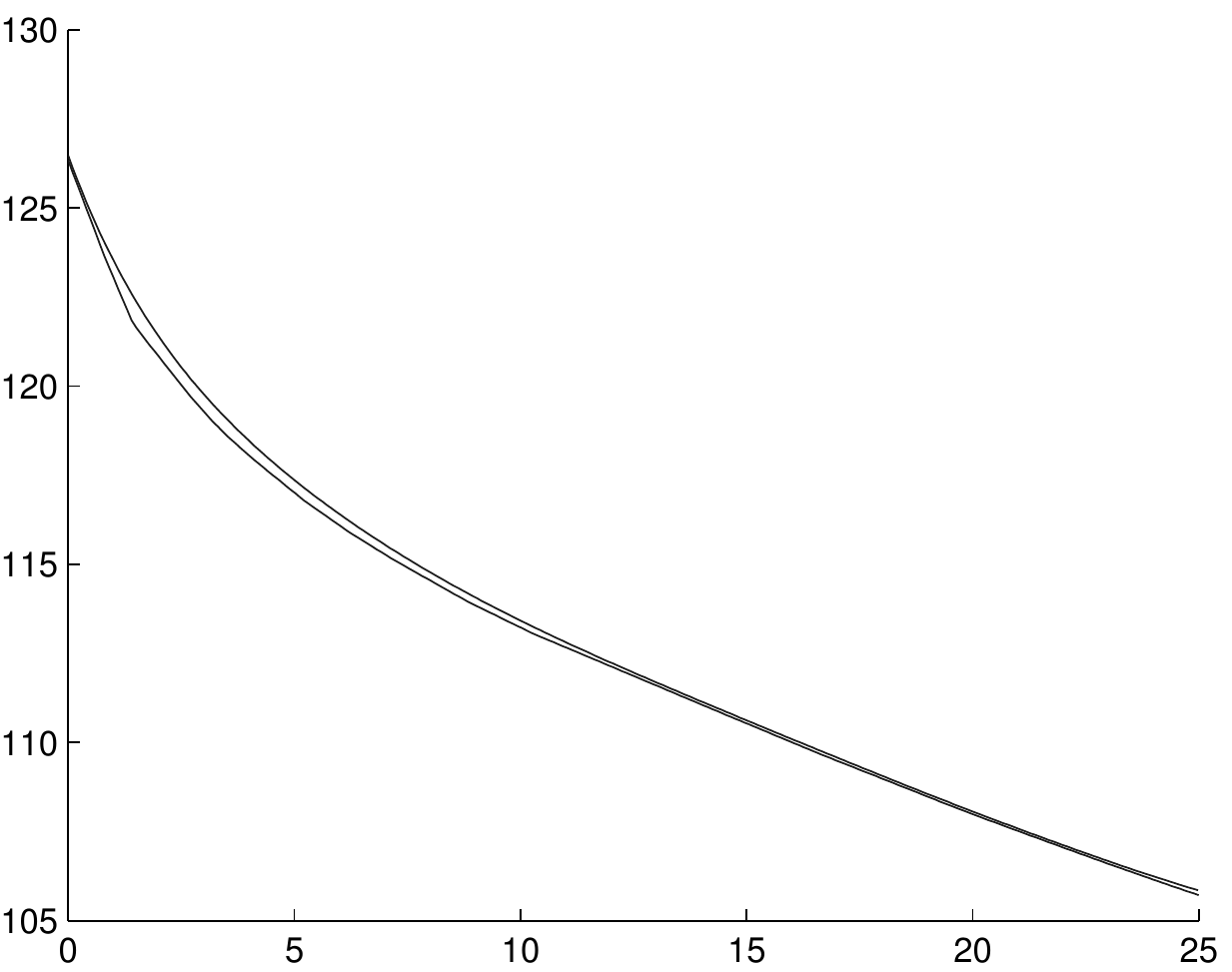}
		\label{fig:ex-dim2-8}
  }
	\subfloat[degree 16]{
		\includegraphics[scale=0.35]{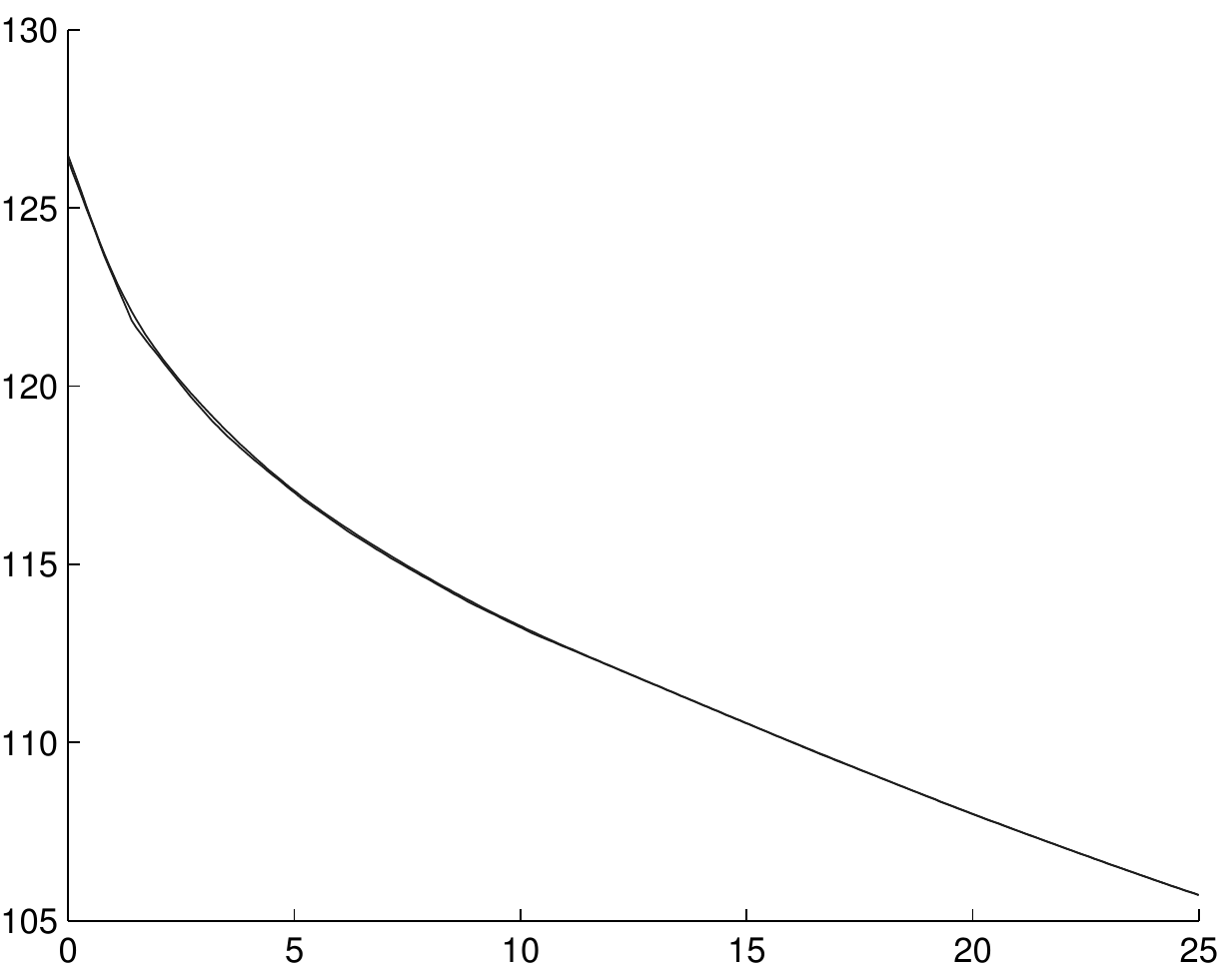}
		\label{fig:ex-dim2-16}
  }
	\caption{Numerical example with two objectives indicating the quality of the inner approximation with a polynomial decision rule. The lowest curve represents the set $\PS$ which is usually not known.}
	\label{fig:ex-dim2}
\end{figure}

We solve linear programs with Matlab linprog. We enter linear constraints with LMI uncertainty into YALMIP \cite{Lofberg:2010}. YALMIP reformulates this problem as an SDP. In Appendix \ref{sec:sdpreformulation} we show how to do this reformulation by hand in case of a polynomial of degree 3. We let YALMIP export the resulting problem, then we reformulate free variables as the difference of two nonnegative variables using CSDP's convertf, and solve the problems with SDPA \cite{SDPA} (SDPA-DD \cite{SDPADD} for the problem with a polynomial of degree 16).
Figure \ref{fig:ex-dim2} shows the resulting solutions. The solution time ranges from 21 seconds for the polynomial approximation of degree 4 to 4 minutes for degree 8 (with SDPA), and 45 minutes for degree 16 (with SDPA-DD).

\subsection{Three objectives} \label{sec:threeobj}
We semi-randomly construct vectors $\vec{c^1} \in [0,1]^{10}$, $\vec{c^2} \in [0,1]^{10}$ and $\vec{c^3} \in [-1,0]^{10}$, and take $\R_+^{10}$ as the feasible region. Recall from Table \ref{tbl:tractability} that we have an exact result for a quadratic inner approximation over an ellipsoidal set, and a conservative result for polynomial inner approximations of arbitrary degree over semialgebraic sets. Again we use YALMIP and SDPA to formulate the problem and solve the resulting SDPs. We use YALMIPs SOS module for constraining an expression to be SOS.

\begin{figure}
	\centering
	\subfloat[degree 2]{
		\includegraphics[scale=0.2]{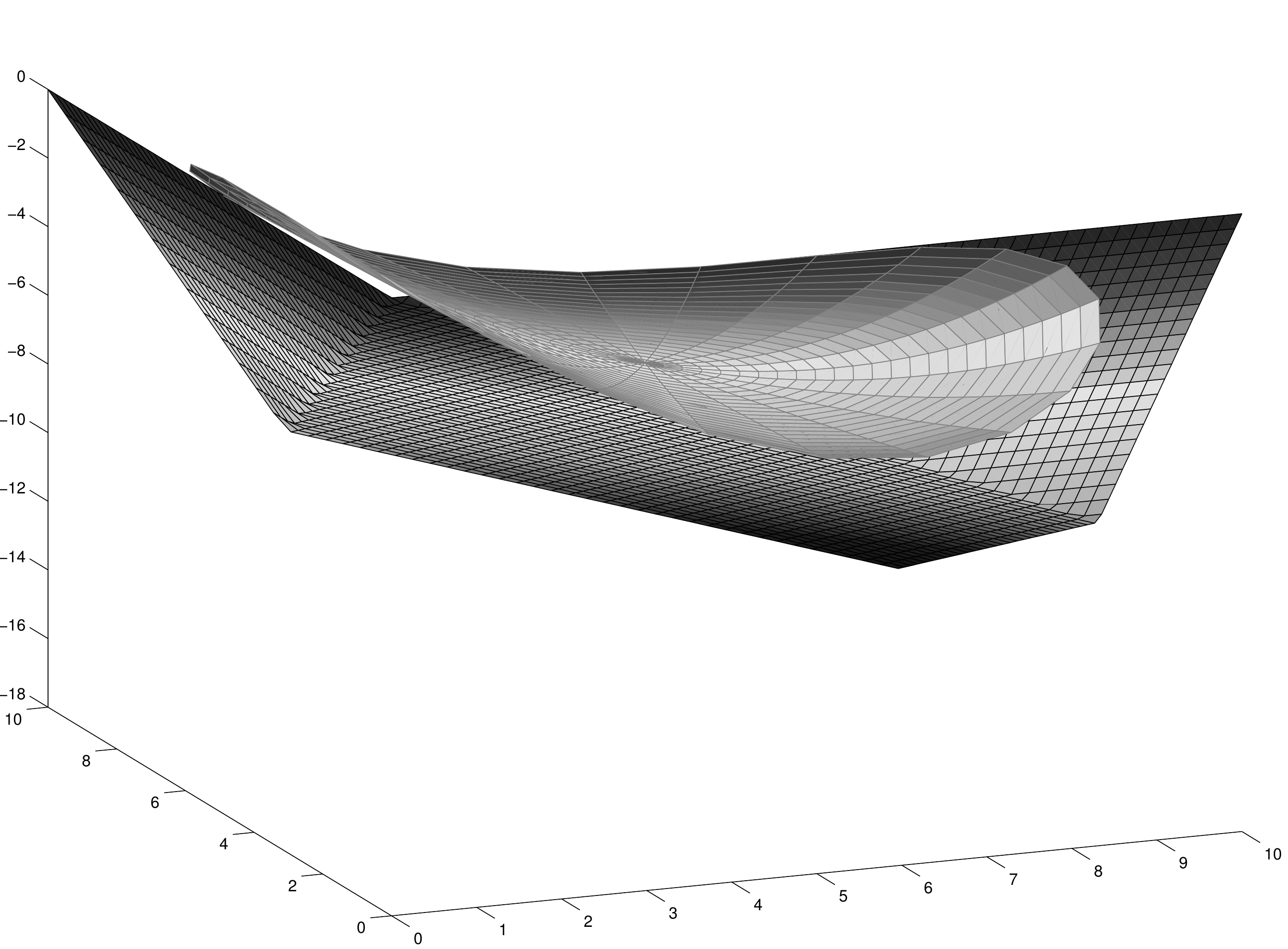}
		\label{fig:ex-dim3-1}
  }
	\quad
	\subfloat[degree 4]{
		\includegraphics[scale=0.2]{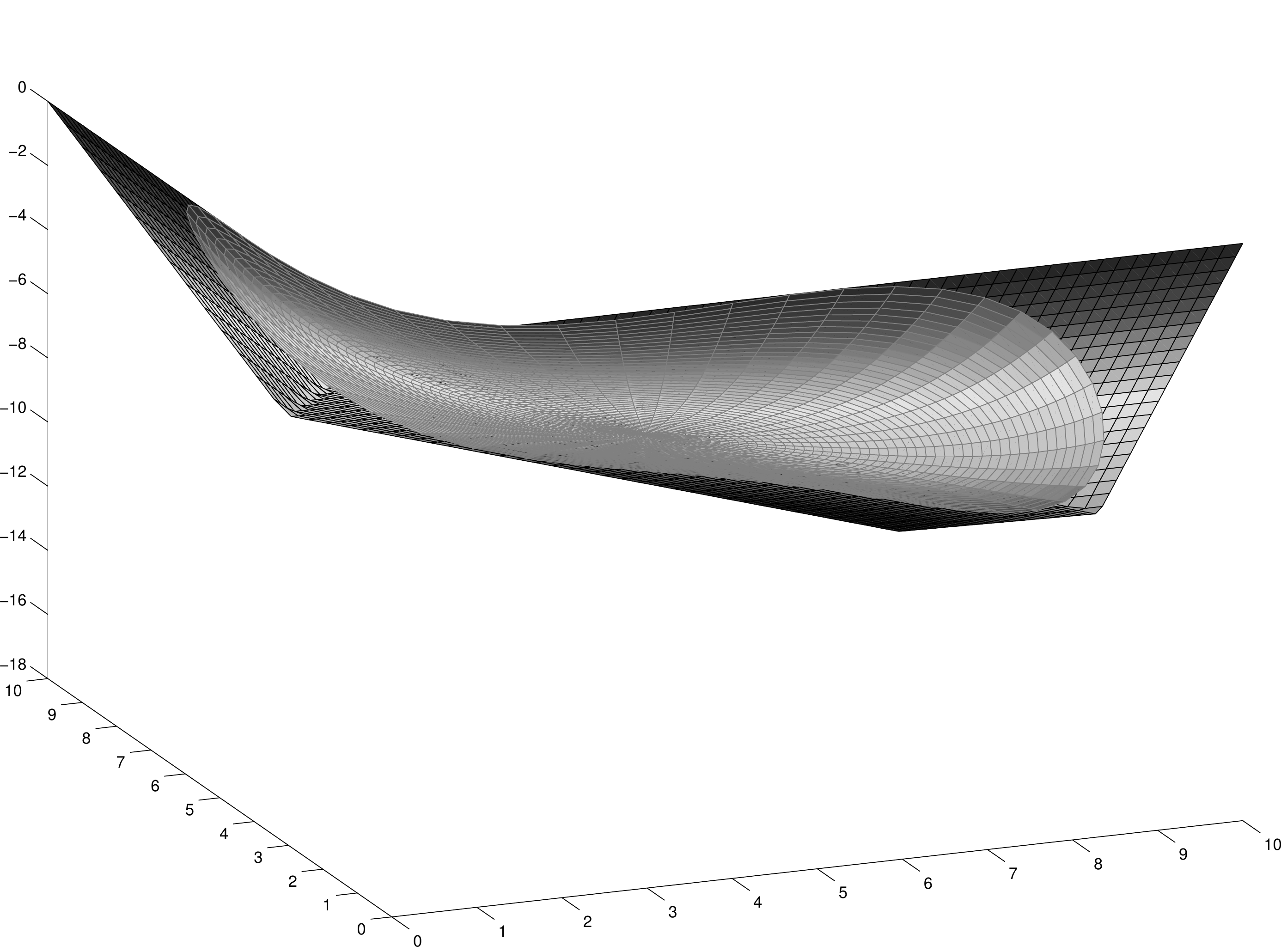}
		\label{fig:ex-dim3-2}
	}
	\caption{Numerical example with three objectives indicating the quality of the inner approximation with a polynomial decision rule. The lowest set represents the set $\PS$ which is usually not known.}
	\label{fig:ex-dim3}
\end{figure}

We take $\{\vec{u} : \twonorm{\vec{u}-\vec{5}} \leq 5\}$ as the area of interest for $(\vec{c^1}) \transp \vec{x}$ and $(\vec{c^2}) \transp \vec{x}$, and approximate the Pareto set with a polynomial of degree 2 and with a polynomial of degree 4. For degree 2, we solve the exact robust counterpart, while for degree 4 we solve the SOS approximation. We also solve the SOS approximation for degree 2, and notice that the inner approximation is the same as with the exact robust counterpart. Figure \ref{fig:ex-dim3} shows that the polynomial of degree 4 gives an approximation that is closer to $\PS$ than the polynomial of degree 2. The solution time is around 1.6 seconds for all three approximations.

\appendix
\section{Outer approximation\label{sec:outerapprox}}
In this appendix we show how ARO can be used to construct an outer approximation. While theoretically possible, it is practically not tractable to determine polynomial approximations of degree 2 or higher due to computational issues. For linear approximations, the results are trivial. We still mention the result because it also uses ARO to approximate $\PS$.

Let $U=[a^1,b^1] \times [a^2,b^2] \times \cdots \times [a^{k-1},b^{k-1}]$ with $[a^i,b^i]$ be the domain of interest for $(c^i) \transp \vec{x}$ ($i=1,2,\ldots,k-1$). We construct the set $\OPS$ by creating a function function $\ell : \R^{k-1} \to \R$ for which $(u^1, u^2, \ldots, u^{k-1}, \ell(u^1, u^2, \ldots, u^{k-1}))$ is in $\PS+\R^k_-$, and optimizing over this function:
\begin{subequations}
\begin{align}
 \max_{\ell} ~ & \int_U \ell(u_1,u_2,\ldots,u_{k-1}) du \label{oa-obj} \\
 \mbox{s.t.} ~ & \ell((\vec{c^1}) \transp \vec{x},(\vec{c^2}) \transp \vec{x},\ldots,(c^{k-1}) \transp \vec{x}) \leq (\vec{c^k}) \transp \vec{x} \notag \\
               & \qquad \qquad  \forall \vec{x} : A\vec{x} \leq \vec{b}, a^i \leq (\vec{c^i}) \transp \vec{x} \leq b^i. \label{oa-feas}
\end{align}
\end{subequations}
Given a solution to this optimization problem, the outer approximation is $\newline \{ (u_1, u_2, \ldots, u_{k-1}, \ell(u_1, u_2, \ldots, u_{k-1})) : \vec{u} \in U \}$. The objective \eqref{oa-obj} maximizes the volume under this approximation. Constraint \eqref{oa-feas} ensures that the outer approximation as a function of $(\vec{c^1}) \transp \vec{x}$ lies under $(\vec{c^2}) \transp \vec{x}$ for every $\vec{x}$ in the domain of interest.

An optimal outer approximation is tangent to the Pareto curve at (at least) one point. This becomes clear from \eqref{oa-feas}: this constraint holds with equality for at least one $x$ because otherwise we can add a constant to $\ell$ without losing feasibility, which contradicts optimality. Previous results force the decision maker to specify either this point of tangency or the derivative at this point a priori. Our formulation determines the point of tangency in such a way that the volume enclosed between this linear outer approximation and the Pareto curve over the set $U$, but in the linear case this turns out to give a trivial result.

When $\ell$ is linear, the problem \eqref{oa-feas} can be reformulated as an LP using ARO. For the case $k=2$ (two objectives) it can be shown that the optimal linear $\ell$ is a line tangent to $\PS$ at $\frac{a^1+b^1}{2}$, i.e. halfway the interval of interest. We conjecture that in higher dimensions the point of tangency is the barycenter of $U$. This would imply that the formulation for the outer approximation is not interesting because it is already known how to obtain an outer approximation that is tangent at a given point.

For nonlinear $\ell$ the SOS framework used in Section \ref{sec:ia} can be used to reformulate the problem as an SDP when $\ell$ is a polynomial of arbitrary degree $d$. This is a polynomial in the vector $\vec{x}$, so the number of terms is ${n+d \choose d}$, which is also the order of the matrix in the SDP. Even for a quadratic function, the size of the SDP is often too large to solve.

The user has to specify the domain of interest $[u^1,u^2]$. Specifying the wrong domain does not lead to infeasibility. However, if the interval is too large, i.e. $(\vec{c^1}) \transp \vec{x}$ does not range through the full interval, part of the outer approximation is meaningless because $\PS$ is inexistent for some $\vec{u}$.
\section{Derivation of the SDP formulation for a polynomial inner approximation with two objectives (online supplement) \label{sec:sdpreformulation}}
We give a derivation of the SDP formulation of \eqref{ia-obj}-\eqref{ia-feas} in case $\vec{x(u)} = \vec{\alpha_0} + \vec{\alpha_1} u + \vec{\alpha_2} u^2 + \vec{\alpha_3} u^3$ (where $\vec{\alpha_i}$ in $\R^{170}$, i=1,2,3) for the numerical example of Section \ref{sec:twoobj}. Suppose the area of interest for $(\vec{c^1}) \transp \vec{x}$ is $[0,25]$, then $u$ runs from 0 to 25. Because in the result by \cite[Lemma 14.3.4]{BenTal} $u$ runs from -1 to 1, we use the following linear transformation $\vec{x} \to D\vec{x}+\vec{d}$ to transform $\Zset = \{(u,u^2,u^3) : -1 \leq u \leq 1\}$ into $\{(u,u^2,u^3) : 0 \leq u \leq 25\}$:
\begin{align*}
D = \begin{pmatrix} 12.5 & 0 & 0 \\ 312.5 & 156.25 & 0 \\ 5859.375 & 5859.375 & 1953.125 \end{pmatrix}
&&
\vec{d} = \begin{pmatrix} 12.5 \\ 156.25 \\ 1953.125 \end{pmatrix}
\end{align*}
The problem can now be written as follows:
\begin{subequations}
\begin{align}
 \min \qquad & (\vec{c^2}) \transp (25 \vec{\alpha_0} + \frac{25^2}{2} \vec{\alpha_1} + \frac{25^3}{3} \vec{\alpha_2} + \frac{25^4}{4} \vec{\alpha_3}) \notag \\
             & (\vec{c^1}) \transp \left(\vec{\alpha_0} + [\vec{\alpha_1} \quad \vec{\alpha_2} \quad \vec{\alpha_3}] (D \vec{\zeta} + \vec{d})\right) \leq  (D \vec{\zeta} + \vec{d})_1  \notag \\
             & \qquad \forall \vec{\zeta} \in \Zset \label{ia-ap-moving} \\
             & A\left(\vec{\alpha_0} + [\vec{\alpha_1} \quad \vec{\alpha_2} \quad \vec{\alpha_3}] (D \vec{\zeta} + \vec{d})\right) \leq \vec{b} \notag \\
             & \qquad \forall \vec{\zeta} \in \Zset, \label{ia-ap-feas}
\end{align}
\end{subequations}
where:
\begin{align*}
\Zset = \{ \vec{\zeta} \in \R^3 : \begin{pmatrix} 1 \\ \vec{\zeta} \end{pmatrix} = \begin{pmatrix} 1 & 0 & 0 & 0 \\ 0 & 2 & 0 & 0 \\ 3 & 0 & 4 & 0 \\ 0 & 4 & 0 & 8 \\ 3 & 0 & 4 & 0 \\ 0 & 2 & 0 & 0 \\ 1 & 0 & 0 & 0 \end{pmatrix} \transp \begin{pmatrix} \lambda_0 \\ \lambda_1 \\ \lambda_2 \\ \lambda_3 \\ \lambda_4 \\ \lambda_5 \\ \lambda_6 \end{pmatrix},  \quad
\begin{pmatrix} \lambda_0 & \lambda_1 & \lambda_2 & \lambda_3 \\ \lambda_1 & \lambda_2 & \lambda_3 & \lambda_4 \\ \lambda_2 & \lambda_3 & \lambda_4 & \lambda_5 \\ \lambda_3 & \lambda_4 & \lambda_5 & \lambda_6 \end{pmatrix} \succeq 0 \}.
\end{align*}

Constraints \eqref{ia-ap-moving} and \eqref{ia-ap-feas} are a total of 151 semi-infinite constraints that have to hold for an infinite number of $\vec{\zeta}$. Let $\vec{A_j}$ denote the $j^{th}$ row of $A$. In order to allow for shorter notation, we define the linear function $\ell : \R^{170} \times \R^{170} \times \R^{170} \to \R^3$ vector $\vec{\ell(\alpha)} := A_j [\vec{\alpha_1} \quad \vec{\alpha_2} \quad \vec{\alpha_3}] D$. The $j^{th}$ constraint of \eqref{ia-ap-feas} can be rearranged to:
\begin{align*}
  \vec{A_j} \vec{\alpha_0} - \vec{b} + [\vec{\alpha_1} \quad \vec{\alpha_2} \quad \vec{\alpha_3}] \vec{d}
+ \vec{\ell(\alpha)} \vec{\zeta} \leq 0  && \forall \vec{\zeta} \in \Zset,
\end{align*}
which is equivalent to:
\begin{align}
  \vec{A_j} \vec{\alpha_0} - \vec{b} + [\vec{\alpha_1} \quad \vec{\alpha_2} \quad \vec{\alpha_3}] \vec{d} + \max_{\vec{\zeta} \in \Zset} \left\{ \vec{\ell(\alpha)} \vec{\zeta} \right\} \leq 0. \label{ex-sdp}
\end{align}
The maximization problem is an SDP. Replacing this problem with its SDP dual and omitting the $\min$ operator is a well-known method to transform a semi-infinite constraint into a single constraint. We show how to do this. In the following SDP, we take the $4 \times 4$ matrix with $\lambda_i$'s in the description of $\Zset$ as our variable $X$. Let $\ea{i}$ denote the $i^{th}$ component of $\vec{\ell(\alpha)}$.

The optimization problem in \eqref{ex-sdp} is:
\begin{align*}
\max          \qquad & \langle C,X \rangle \\
\textrm{s.t.} \qquad & \langle A_i, X \rangle = b_i \quad (i=1,2,3,4) \\
                     & X \succeq 0,
\end{align*}
where $\langle \cdot,\cdot \rangle$ denotes the trace inner product, $b_1=1$, $b_2=b_3=b_4=0$, and:
\begin{align*}
C = \begin{pmatrix}
      \ea{3}       & \ea{2}     & g(\alpha)  & \ea{2}    \\
      \ea{2}       & g(\alpha)  & \ea{2}     & g(\alpha) \\
      g(\alpha)    & \ea{2}     & g(\alpha)  & \ea{2}    \\
      \ea{2}       & g(\alpha)  & \ea{2}     & \ea{3}
    \end{pmatrix},
\end{align*}\begin{align*}
\mbox{with} ~ g(\alpha) = \frac{4}{3}\ea{1}+\ea{3} \mbox{, and:} \\
A_1 = \begin{pmatrix} 0 & 0 & 0 & 2 \\ 0 & 0 & 2 & 0 \\ 0 & 2 & 0 & 0 \\ 2 & 0 & 0 & 0\end{pmatrix}& , &
A_2 = \begin{pmatrix} 0 & 0 & -\frac{1}{2} & 0 \\ 0 & 1 & 0 & 0 \\ -\frac{1}{2} & 0 & 0 & 0 \\ 0 & 0 & 0 & 0\end{pmatrix}, \\
A_3 = \begin{pmatrix} 0 & 0 & 0 & -1 \\ 0 & 0 & 1 & 0 \\ 0 & 1 & 0 & 0 \\ -1 & 0 & 0 & 0\end{pmatrix} & , &
A_4 = \begin{pmatrix} 0 & 0 & 0 & 0 \\ 0 & 0 & 0 & -\frac{1}{2} \\ 0 & 0 & 1 & 0 \\ 0 & -\frac{1}{2} & 0 & 0\end{pmatrix}.
\end{align*}
By formulating the dual and putting this into constraint \eqref{ex-sdp}, we get the following robust counterpart:
\begin{align*}
\vec{A_j} \vec{\alpha_0} - \vec{b} + [\vec{\alpha_1} \quad \vec{\alpha_2} \quad \vec{\alpha_3}] \vec{d} + y_1 \leq 0 \\
\sum_{i=1}^4 y_i A_i - C \succeq 0.
\end{align*}
Constraint \eqref{ia-ap-moving} and the other constraints \eqref{ia-ap-feas} can be transformed in a smilar way.
\section*{Acknowledgments}
We would like to thank R. Sotirov from Tilburg University (The Netherlands) for her input on solving SDPs.
\bibliographystyle{abbrvnatnew}
\bibliography{cleanpareto}
\end{document}